\documentclass[12pt]{article}
\usepackage{amssymb}
\usepackage{amsthm}
\usepackage{epsf,epsfig,amsfonts,graphicx,color}
\usepackage{a4wide}
\newtheorem{theorem}{Theorem}

\newtheorem{remark}{Remark}

\newtheorem{corollary*}{Corollary}
\newtheorem{example}{Example}

\newcommand{\be}{\begin{equation}}
\newcommand{\comp}{\mathbb{C}}
\newcommand{\ee}{\end{equation}}

\newcommand{\R}{\mathbb{R}}
\newcommand{\rn}{\mathbb{R}^n}
\newcommand{\weg}[1]{}
\title{Finsler Conformal Lichnerowicz-Obata conjecture}
\author{Vladimir S. Matveev \and Hans-Bert Rademacher 
\and  Marc Troyanov \and  Abdelghani  Zeghib}
\date{}
\begin{document}
\maketitle
\begin{abstract}
We  prove the Finsler analog of the conformal Lichnerowicz-Obata  
conjecture showing that  a complete and essential conformal vector field on a 
non-Riemannian Finsler manifold is a homothetic 
vector field of a Minkowski metric. 
\\[2ex]
\\[2ex]
\end{abstract}

\section{Definitions and results} 
In this paper a \emph{Finsler metric} on a smooth manifold  $M$ 
is a 
function $F:TM\to \mathbb{R}_{\ge 0} $ satisfying
the following properties:
\begin{enumerate} 
\item It is smooth on $TM\setminus TM_0$, 
where $TM_0$ denotes the zero section of $TM$,
\item For every $x\in M$,  the restriction $F_{|T_xM}$ is a norm on 
$ {T_xM}$, i.e.,  for every $\xi,\,  \eta \in T_xM$ and for every  
nonnegative  $\lambda\in \mathbb{R}_{\ge 0}$   we have 
 \begin{enumerate} 
   \item $F(\lambda \cdot \xi) = \lambda\cdot   F(\xi) $, \label{2a}
   \item $F(\xi+ \eta ) \le F(\xi) + F(\eta)$, \label{2b} 
   \item $F(\xi)= 0 $ $ \Longrightarrow$ $\xi=0$.  \label{2c}
  \end{enumerate}  
  \end{enumerate} 
We do not require  that (the restriction of) 
the function $F$ is strictly convex. In this point our definition 
is more general than the usual definition.
In addition we do not assume the metric to be reversible, i.e.,  
we do not assume that $F(-\xi)= F(\xi).$
Geometrically speaking a Finsler metric is characterized
by a smooth family $ x \in M \mapsto \{ \xi \in T_xM \mid 
F(\xi)=1\}\subset T_xM$ of convex hypersurfaces (sometimes called \emph{indicatrices}, 
cf. \cite{Br}) containing the zero section 
in the tangent bundle.\\

Recent references for Finsler geometry  include 
\cite{BCS,Sh,BBI,Alv}. Particular classes of Finsler metrics
which occur in our results are the following:
\begin{example}[Riemannian metric]  
\label{riem} \rm
For every Riemannian metric $g$ on $M$  
the function $F(x, \xi):= \sqrt{g_{(x)}(\xi, \xi)}$ is a 
Finsler metric. Geometrically the Finsler metric is a smooth
family of ellipsoids.
\end{example} 
\begin{example}[Minkowski  metric]  
\label{mink}
\rm 
Consider a norm  on   $\mathbb{R}^n$, i.e.,  
a function $p:\mathbb{R}^n\to \mathbb{R}_{\ge 0}$ 
satisfying  \ref{2a}, \ref{2b}, \ref{2c}. We canonically 
identify $T\mathbb{R}^n$ with $\mathbb{R}^n\times  \mathbb{R}^n$ 
with coordinates 
$\bigl(\underbrace{(x_1,...,x_n)}_{x\in \mathbb{R}^n}, 
\underbrace{(\xi_1,...,\xi_n)}_{\xi\in T_x\mathbb{R}^n}\bigr)$. 
Then, $F(x, \xi) := p(\xi)$ is a Finsler metric.    
The metric is translation invariant, 
on the other hand every  translation invariant 
Finsler metric is a Minkowski metric. 
Due to the translation invariance the Finsler metric
is uniquely characterized by a convex hypersurface 
in a single tangent space $T_xM$. 
\end{example}
Two Finsler metrics $F$ and $F_1$ on
an open subset $U\subseteq M$ are called 
\emph{conformally equivalent,} if $F_1= \lambda\cdot F$ 
for a  nowhere vanishing function $\lambda  $ on $U$. 
We say that a  differentiable mapping
$f: (M_1, F_1) \to (M_2, F_2)  $ 
is \emph{conformal}, if  the pullback of the metric $F_2$ 
is conformally equivalent to $F_1$, i.e., 
if for every $\xi \in T_xM$ we have
$F_2\left(df_x(\xi)\right)=\lambda (x) F_1\left(\xi\right).$
If the conformal factor $\lambda$ is constant the map is 
called {\em homothetic,} for $\lambda=1$ it is {\em isometric.}
A vector  field is called \emph{conformal} (resp. 
{\em homothetic} or {\em isometric}) 
if its local flow acts by \emph{conformal} (resp.
homothetic or isometric)
local diffeomorphisms. If the conformal
vector field $v$ is 
{\em complete} then the flow 
$\phi^t: M \rightarrow M, t\in \R$
of $v$ is a one-parameter
group of conformal
diffeomorphisms of the manifold $M.$
\\
Obviously, if a metric $F_1$ is conformally equivalent to $F$, 
then every conformal vector field for $F$ is also a conformal 
vector field for $F_1$. 
\begin{example} 
\label{exa:eucl}
\rm
For the Finsler metric $F:= \sqrt{g(\xi, \xi)}$ 
from Example \ref{riem},  conformal vector fields for the   
Riemannian metric  $g$ 
are  conformal vector fields  for the Finsler metric $F$, 
and vice versa. For Euclidean space $\R^n$ the description of 
conformal mappings 
for $n=3$ is due to Liouville \cite{Li} and for $n \ge 3$ to Lie
\cite{Lie}, for recent expositions cf. for
example \cite[Thm. A.3.7]{BP}, \cite{Ku} and \cite{KR}.
\end{example}
\begin{example} 
\rm
For the Minkowski Finsler metric 
$F(x, \xi) := p(\xi)$ from Example \ref{mink}, the 
mappings $x \in \R^n
\mapsto H_t(x)=t\cdot x\in \R^n$ are homotheties
for all $t>0.$ Then,  the vector field
$v(x)= \left.\frac{d}{dt}\right|_{t=0}H_t(x)
 = \sum_{i=1}^n x_i\frac{\partial }{\partial x_i} $  
is the corresponding homothetic vector field.
\end{example} 
Now we can state our main result:
\begin{theorem}
\label{thm:main}
Suppose $v$ is a conformal 
and complete vector field on a connected Finsler 
manifold $(M, F)$ of dimension $n \ge 2.$
Then, at least one of the following statement  holds. 
\begin{enumerate} 
\item There exists a Finsler  metric $F_1$ 
conformally equivalent to  $F$ 
such that the flow of $v$ preserves the Finsler metric $F_1.$ 
\item 
The manifold 
$M$ is conformally equivalent to the sphere $S^n$
with its standard Riemannian metric.
\item The manifold $M$ is  diffeomorphic to 
Euclidean space $\mathbb{R}^n$,  
and the Finsler metric $F$ is conformally equivalent 
to a Minkowski metric, cf. Example \ref{mink}.
The vector field $v$ with respect to the Minkowski metric
is homothetic.
\end{enumerate} 
\end{theorem}
%
For Riemannian metrics, the statement above is called  the 
conformal \emph{Lichnerowicz-Obata} conjecture, 
and was proved independently by D.Alekseevksii~\cite{Al},
J.Ferrand~\cite{Fe2}, M.Obata~\cite{Ob}  and R.Schoen~\cite{schoen}, 
see also \cite{La,FT}.  
Of course, in the Riemannian case, 
Example \ref{mink} corresponds to the Euclidean  
metric on $\mathbb{R}^n$.
A conformal vector field satisfying
the assumptions of the first case is also called {\em inessential,}
otherwise it is called {\em essential.}
\begin{remark}\rm
\label{rem:essential}
Theorem~\ref{thm:main} also implies the following result:
If the conformal group is essential, i.e. if there is no
conformally equivalent metric such that the conformal group
becomes the isometry group, 
then the metric is conformally equivalent to the round sphere,
or to a Minkowski space.
\end{remark}
Theorem \ref{thm:main} was announced  
in \cite{Ze} under the following  additional assumptions: 
$M$ is closed, and the Finsler metric $F$ is \emph{strictly convex}, 
i.e., the second derivative of $F^2_{|T_pM}$ 
has rank $n-1$  at every point on $T_pM -TM_0.$  
The proof is sketched in \cite{Ze}. 
It  is long  and actually is a repeating of the proof from 
\cite{Fe2} (which is technically very nontrivial)  
in the Finsler case.

Our proof of Theorem \ref{thm:main} is much shorter. 
It is  based on the following observation: 
for every Finsler metric $F$ we can canonically  
construct a Riemannian metric $g$ such that if $v$  
is a conformal vector field for $F$, 
then it also a  conformal vector field  for $g$. 
Then, by the Riemannian version of Theorem \ref{thm:main},  
the following two cases are possible: 
\begin{itemize} 
\item  The flow of $v$ acts by  isometries of a certain 
Riemannian metric $g_1$ conformally equivalent to the  
Riemannian metrics $g$.  This case will be called  
``trivial case'' in the proof of Theorem~\ref{thm:main}. 
In this case,  it immediately follows, that the  flow of  
$v$ acts by the isometries of a particular  metric $F_1$ 
conformally equivalent to $F$. 
\item The manifold  is $S^n$ or $\mathbb{R}^n$, and the metric $g$ 
is conformally equivalent to the standard 
metric. In this case, all possible essential conformal  
vector fields  $v$  can be explicitly described,
cf. Example~\ref{exa:eucl}.
A direct analysis of the flow of such vector field shows, 
that the only Finsler metrics  for which $v$ is a conformal 
vector field are as in Theorem~\ref{thm:main}.
\end{itemize}
\begin{remark}
\label{rem:surface} 
\rm
In conformal geometry the case of surfaces 
$n=2$ is
special due to the existence of holomorphic functions. 
Any holomorphic function defined
on an open subset
in the complex plane $\mathbb{C}$ with everywhere
non-vanishing derivative is conformal.
This shows that the part of 
Liouville's theorem on conformal transformations
of Euclidean spaces stating
that a conformal diffeomorphism between open subsets
of Euclidean space is the restriction of a conformal
diffeomorphism of the standard sphere only holds
for dimensions $n\ge 3.$ On the other hand the
description of the conformal diffeomorphism of
the $n$-dimensional sphere $S^n$ as compositions of 
homotheties and inversions in the Euclidean
space $\R^n\cong S^n-\{p\}$ also holds for $n=2,$
as one concludes from the standard classification of
(anti)holomorphic functions on $\comp$ resp.
$\comp P^1.$
It is shown by Alekseevskii~\cite[Thm.8]{Al}
that an essential and complete conformal vector field
on a surface only exists on the $2$-sphere with
the standard metric or on Euclidean $2$-space.
Therefore for our main result the case $n=2$
is not exceptional.
\end{remark}
\section{Averaged Riemannian metric}  \label{sec2} 
For a given smooth norm $p$  on $\mathbb{R}^n$ we construct canonically 
 a positive definite  symmetric  bilinear form 
$g:\mathbb{R}^n\times  \mathbb{R}^n \to \mathbb{R}$. 
 
 For a Finsler metric $F$, the role of $p$ will play the restriction of $F$ to $T_xM$.  
 We will see that the constructed  $g$ will smoothly depend on $x$, i.e., 
$g_{(x)}$ is a Riemannian metric.

Consider the  unit sphere $S_1 =  \{ \xi \in \mathbb{R}^n \mid p(\xi)=1\}$
of the norm $p.$  
Consider the (unique)  volume form $\Omega $ on $\mathbb{R}^n$  
such that the volume of the 1-ball  $B_1 =  \{ \xi \in \rn  \mid p(\xi)\le 1\}$ 
equals $1$.
 
 Denote by $\omega$ the volume form  on $S_1$,   whose values on the  vectors 
$\eta_1,...,\eta_{n-1}  $ tangent to $S_1$ at the point $\xi\in S_1$ are  
 given by  
$\omega(\eta_1,...,\eta_{n-1}) :=   \Omega(\xi,\eta_1,\eta_2,...,\eta_{n-1})$.

Now, for every point $\xi \in S_1$, consider the symmetric  bilinear 
form $b_{(\xi)}: \mathbb{R}^n  \times \mathbb{R}^n   \to 
\mathbb{R} $, $b_{(\xi)} (\eta, \nu) =   D_{(\xi)}^2p^2 (\eta, \nu)$.
In this formula, $D_{(\xi)}^2p^2$ is the second differential at the point 
$\xi$ of the
function $p^2$ on $\mathbb{R}^n$.  The analytic expression for $b_{(\xi)} $ 
in the coordinates  $(\xi_1,..., \xi_n)$
is 
  \begin{equation}  b_{(\xi)} (\eta, \nu) = \sum_{i,j} 
\frac{\partial^2  p^2 (\xi)}{\partial \xi_j\partial \xi_j} \eta_i \nu_j. \label{p2} 
  \end{equation}       
Since the norm $p$ is  convex, the bilinear form (\ref{p2})  
is nonnegative definite: for all $\eta$ we have  
\begin{equation}  b_{(\xi)} (\eta, \eta) \ge 0. \label{p12} 
  \end{equation}  
   Clearly, for every $\xi\in S_1 $, we have  \begin{equation}   
b_{(\xi)} (\xi, \xi)> 0\label{p4} \end{equation}

Now consider  the following  bilinear symmetric $2-$form 
$g$ on $\mathbb{R}^n$:  for $\eta, \nu \in \rn$,  we put 
$$g(\eta, \nu) = 
 \int_{S_1} b_{(\xi)}(\eta, \nu) \omega.
  $$
We assume that the   orientation of $S_1$  is chosen in such a way that 
$ \int_{S_1} \omega\ge 0$. Because of  (\ref{p12})  and    (\ref{p4}), 
$g$  is positive definite.  

\begin{remark} \label{last}
\rm
If the norm $p$ comes from a scalar product, i.e., if 
$p(\xi)=  \sqrt{b_1(\xi,\xi)}$ for a certain positive definite symmetric 2-form
$b_1$, then $b$  is equal to $b_1$ multiplied by a  constant only depending
on the dimension.
\end{remark} 

Starting with a Finsler metric $F$,  we can use this construction
for every tangent space $T_xM$ of the manifold, the role of $p$ is played
by the restriction $F_{|T_xM}$ of the Finsler metric to the tangent space
$T_xM.$ 
Since this construction depends smoothly on the point  $x\in M$, we obtain
a Riemannian metric $g=g(F)$ on $M.$ We call this metric the \emph{averaged
Riemannian metric} of the Finsler metric $F.$
   
\begin{remark} 
\rm
\label{lemma1} \label{lem:eins}
    It is easy to check that for   the metric $F_1:= \lambda(x) \cdot F$  
 the constructed  metric $g_1$ is conformally equivalent to the    
metric $g$ constructed for $F$. More precisely,  $g_1= \lambda(x)^2 \cdot g$.  
Then, a conformal diffeomorphism 
 (conformal vector field,  respectively)  for 
  $F$ is also a conformal diffeomorphism  (conformal vector field, respectively) 
for  $g$. Moreover, if $v$ is conformal for $F$ and is an 
isometry (homothety, respectively) for $g$,
   then it is an isometry (homothety, respectively) for $F$ as well. 
\end{remark} 
\begin{remark}
 \label{rem:szabo}
\rm
This averaging construction is quite  natural   and it is very possible  
that other researchers in Finsler geometry already thought
about it, but we could not find any reference about it in the  
literature, nor any significant result in Finsler geometry
whose proof is based on the averaged metric. It would certainly be  
worthwile to further investigate its properties.
Recently, Szabo~\cite{Sz} uses a similar averaging construction
to explicitely construct all Finlser  Berwald metrics. 
There are other canonical constructions of a bilinear form  starting from a norm. 
R. Schneider told us, that a for a  convex geometer  the  natural bilinear form 
corresponding to a convex body is one corresponding to the John ellipsoid of this convex body. 
\\
These constructions  have the nice  properties listed in Remarks \ref{last}, \ref{lem:eins}. 
We still prefer our averaged Riemannian metric,  since the method of Szabo assumes that 
the norm is strictly convex, and  since it is not clear whether the John ellipsoid  depends 
smoothly on the norm.
In~\cite{To} Torrome suggests another averaging construction for
Finsler metrics.
\end{remark}

\section{Proof  of Theorem 1}

Let $v$ be a  complete conformal vector field on a 
connected  Finsler manifold $(M, F)$. 
Then, it is also a conformal vector field for the 
averaged Riemannian metric $g$. Then, by the Riemannian 
version of our Theorem, which, as we explained in 
the introduction, was proved in \cite{Ob,Al,Yo,Fe2,schoen},  
we have  the following possibilities: 

\begin{itemize} 
\item[] {(\bf Trivial case}) $v$ is a Killing vector field 
of a conformally equivalent metric $\lambda(x)^2 g$. 

\item[] {(\bf Interesting  case}) For a certain function 
$\lambda$, 
the Riemannian manifold  $(M^n, \lambda(x)^2  g)$ 
is $(\mathbb{R}^n, g_{0})$, or $({S}^n, g_{{1}})$, 
where $g_0$ resp. $g_{{1}}$ is the 
Euclidean metric on $\rn$ resp. the standard  
metric of sectional curvature 1 on $S^n$. \end{itemize}

In the {\bf trivial case}, as   we explained in Remark \ref{lem:eins},  
for a certain function  $\lambda$,  
$v$  is a Killing vector field for the metric 
$F_1:= \lambda  \cdot F$, which was one of the possibilities in 
Theorem \ref{thm:main}. 

Now we treat the {\bf interesting case}. Without loss  of 
generality, we can assume that $(M,g)$  is 
$(\mathbb{R}^n, g_0)$, or $(S^n, g_1).$
\subsection{Case 1: $(M, g)= (\mathbb{R}^n, g_0)$.} 
\label{sub:one}
Since the vector field is complete,
 it generates a
one parameter group $\phi^t: \R^n \rightarrow \R^n$ of 
conformal transformations
with respect to the Finsler metric $F$ and the averaged 
Riemannian metric $g_0.$
It follows from Liouville's theorem that for any $t$
the mapping $\phi^t$ is a homothety of the Riemannian metric $g_0.$
In other words, in an appropriate  cartesian 
coordinate system $(x_1,...,x_n)$,  the conformal diffeomorphism 
$\phi=\phi^1 $ has the form 
\begin{equation} 
\label{phi} 
\phi (x_1,...,x_n)= \mu \cdot (x_1,...,x_n)A,
\end{equation}  
where $A$ is an orthogonal $(n\times n)$-matrix. 
Without loss of generality we can assume that $0<\mu<1$.    
We will show that in this case the metric $F$ is as in 
Example \ref{mink}. 

We identify  
$ T_x \rn $ and $\rn \times \rn$ with the help of the  
cartesian coordinates $x=(x_1,...,x_n)$. 
We assume that  the first component of the product  
$\rn \times \rn$  corresponds to our manifold $\rn$, and that the 
second component of the product  $\rn \times \rn$  corresponds to the tangent spaces. 
The coordinates on the tangent spaces    will be denoted by $\xi$, so
$\bigl(\underbrace{(x_1,...,x_n)}_{x\in \mathbb{R}^n}, 
\underbrace{(\xi_1,...,\xi_n)}_{\xi\in T_x\mathbb{R}^n}\bigr)$ 
is a coordinate system on $ T_x \rn \cong \rn \times \rn$.

Clearly,
 the differential of the mapping $\phi$  given by (\ref{phi}) is given by 
$$d\phi_x(\xi)=\left(\mu\cdot(x_1,...,x_n)  A ,
\mu\cdot(\xi_1,...,\xi_n)  A\right).$$
 Then,  for every $\xi, \eta \in T_x\rn$, we have  
$g_{\phi(x)}( d\phi_x(\xi), d\phi_x(\eta))= \mu^2\cdot g_{(x)}(\xi,\eta)$. 
Hence, by Remark \ref{lem:eins}, 
$F\left(\phi(x), d\phi_x(\xi)\right) = \mu\cdot  F \left(x, \xi\right).$ 
Consider  the mapping $$h: T \rn \cong \rn\times \rn
\rightarrow\rn\times \rn, \  \ \  
h(x_1,...x_n,\xi_1,...,\xi_n)=\left(\mu\cdot  (x_1,...,x_n) A,
(\xi_1,...,\xi_n)A\right).$$ By construction,  this mapping  
satisfies $F\left(h (x,\xi)\right)= F(x,\xi).$
Since the orthogonal group $\mathbb{O}(n)$ is compact,  we can
choose a sequence $m_j \to \infty$ such that
$A^{m_j} \to 1\in \mathbb{O}(n)$ for $j \to \infty.$ Then,   
 $(0,\xi)=\lim_{j\to \infty} h^{m_j}(x,\xi)$. Hence,  
$$F(0,\xi )=F \left(\lim_{j\to \infty} h^{m_j}(x,\xi)\right)= 
\lim_{j\to \infty}F \left( h^{m_j}(x,\xi)\right)=
F  (x,\xi).$$ 
Thus, $F$ is translation invariant and therefore a Minkowski metric, cf.
Example~\ref{mink}. 
Hence in this case,  up to conformal equivalence,  the Finsler metric
is a Minkowski metric,  and the conformal vector field is homothetic.
\subsection{Case 2:
$(M,g) =\left({S}^n,g_1\right)$} 

Then, by \cite[Thm. 12]{La} any
essential
conformal vector field 
$v$ vanishes at exactly one (Case 2a) or  exactly two (Case 2b) points. 
We denote by
$v^{-1}(0)=\{x \in M \mid  v(x)=0\}$ the set of zeroes. 
If we assume $v(x)=0$ we use the stereographic projection
$s_x: S^n-\{x\}\rightarrow \rn$ and obtain with the push forward of
the vector field $v$ a complete and conformal vector field on $\rn.$

\subsubsection{Case 2a: Suppose  $v^{-1}(0)=\{x,y\}, x \not=y$}
\label{3}
Suppose the conformal vector field $v$ vanishes precisely at 
two points $x$ and $y$ of the sphere.  
We will show that the Finsler metric $F$ is in fact Riemannian.  
 
Denote by $s_+:\left(S^n - \{ x\}, g_1\right)\rightarrow (\R^n, g_0)$ 
the stereographic projection from 
$x$ which is conformal with respect to the standard Riemannian
metrics $g_0,g_1$ 
with conformal factor $\sigma_+.$
Here,  $\rn$ 
should be identified with the hyperplane through the origin parallel to
the tangent spaces $T_{ x}S^n$. Then we define a Finsler metric
$F_+$ by $s_+^*\,F_+=\sigma_+ \, F.$ Then the averaged Riemannian
metric of $F_+$ coincides with the Euclidean metric $g_0.$
The push foward vector field $v_+:=s_+^* v$ is a conformal 
and complete vector on $\rn$ with respect to the Finsler metric
$F_+$ as well with respect to the standard metric $g_1.$
This vector field has exactly one zero on $\rn.$
Therefore, by section~\ref{sub:one}, the Finsler metric $F_+$ is 
a Minkowski metric, i.e., translation invariant. In particular we
can assume without loss of generality
that the zero point of $v_+$ is the origin of
$\rn$. Hence we can assume that the zero points of $v$ on $S^n$
are antipodal points, i.e., $v^{-1}(0)=\{\pm x\}.$

The stereographic projection
$s_{-}: \left(S^n - \{ -x\}, g_1\right)\rightarrow (\R^n, g_0); s_-(q)=s_+(-q)$
from $-x$ is a 
conformal mapping with conformal factor $\sigma_{-}$ with
$\sigma_+(-q)=\sigma_-(q), q \in S^n$
i.e., $s_{\pm}^* \,g_0=\sigma_{\pm}^2 \,g_1.$
Then we define also the Finsler metric $F_{-}$ on $\rn$
by 
$s_{-}^* \,F_{-}=\sigma_{-} \,F.$
The averaged Riemannian metrics of $F_{-}$ 
equals the Euclidean metric $g_0.$
The push-forward $v_{-}:=(s_{-})_* v$
of the vector field $v$
is a conformal vector field on $\rn$ with respect to the 
Finsler metric $F_{-}$ and, hence, 
with respect to the standard metric $g_0.$
Both vector fields $v_{\pm}$ are evidently complete 
and have precisely one zero at the origin.
Therefore, by section~\ref{sub:one},
the Finsler metrics $F_{\pm}$ are Minkowski
metrics, i.e., translation invariant.

It is well known that
the composition $s_- \circ s_+^{-1}: \R^n -\{0\} \rightarrow \R^n$
equals the inversion $I(q)=q/g_0(q,q)$ at the unit sphere. Therefore,  the
inversion defines a conformal transformation 
$I: (\R^n-\{0\}, F_+) \rightarrow 
(\R^n-\{0\}, F_-)$ between  the two Minkowski metrics.
The differential $d I_q$ of the inversion 
at a point $q \in S^{n-1}:=\{u \in \R^n \mid 
g_0(u,u)=1\}$ equals the reflection $R_q$ at the hyperplane normal to $q.$ 
This implies that  $d I_q^* F_+=R_q^* F_+=F_-$ for any $q \in S^{n-1}.$ 
Since the reflections
generate the orthogonal group and since the Finsler metrics $F_{\pm}$ are
translation invariant, it follows that the norms $F_{\pm}|T_0M$ 
at the origin are invariant under the full orthogonal group and hence 
Euclidean. 

\subsubsection{Case 2b:\enspace $v^{-1}(0)=\{x\}$}
We assume that the vector field $v$ on $S^n$
vanishes precisely at one  point $x\in S^n$. 
We will  again  show that the metric $F$ is  Riemannian.     

We again consider the stereographic 
projections 
$s_{\pm}:S^n-\{\pm x\}\rightarrow \rn$ 
from the points
$x, -x$ as introduced in Section~\ref{3},  and
denote by
$F_{\pm}:=\left(s_{\pm}\right)_*F$ the
induced Finsler metrics on $\rn.$
The push-forward $v_+$ of $v$ with respect to
$s_+$ vanishes nowhere on
$\rn$ and is complete, let $\psi^t$ be its flow 
on $\R^n.$ Then Liouville's theorem implies that
for an arbitrary $t$ the conformal diffeomorphism
$f=\psi^t$ has the form $f(x)=\mu A x+b$
with an orthogonal matrix $A$ and 
$\mu >0, b \in \R^n.$ Since the mapping $f$ has no fixed point
it follows that $b \not=0\,;\, \mu=1$ and $A b =b.$
We introduce the following notation:
$f_{A,b}(q)=Aq+b$ for an orthogonal matrix $A$ and
$b \in \rn$ with $Ab=b.$ 

If we use the stereographic projection
$s_-$,  then the push-forward of $v$ has a zero 
in the origin $0$ and the mapping
$f$ transforms to
$\overline{f_{A,b}}= I \circ f_{A,b} \circ I$ where 
$I= \sigma_-\circ \sigma_+^{-1}$ 
is the inversion at the unit sphere.
Hence $\overline{f}_{A,b}(q)=
\frac{A q + b\|q\|^2}{1+2 \left<Aq,b\right>
+\|b\|^2\|q\|^2}$ where $<.,.>=g_0(.,.)$ with
related norm $\|.\|.$ The conformal factor is
given by $\psi(q)=\frac{1}{1+2 \left<Aq,b\right>
+\|b\|^2\|q\|^2}.$ In particular the conformal mapping
$\overline{f}_{A,b}$ induces at the fixed point
$0$ the map 
\begin{equation}
\label{eq:fnull}
\xi \in T_0\rn \mapsto 
d\,\left(\overline{f}_{A,b}\right)_0(\xi)=
A \xi \in T_0\rn
\end{equation}
which is an isometry also with respect
to the restriction of the  Finsler metric $F_-$ to $0$ 
since $\psi(0)=1.$

For an orthogonal mapping $A$ we introduce the
map $h_A: z \in \rn \rightarrow A z \in \rn$
with induced mapping $(z,\xi)\in T_z\rn=\rn\times \rn
 \mapsto dh_A (z,\xi)=(Az,A\xi) \in T_z\rn\,.$
We want to show that the map $h_A$ is an isometry for
the Finsler metric $F_-.$ 
Let $v_1$ be the vector field on $S^n$ which 
corresponds to the parallel vector field $b$
on $\rn$ with respect to the stereographic projection
$s_+,$ i.e., $ds_+(v_1)(q)=b$ for all $q \in S^n.$
The vector field $v_1$
is a conformal vector field with respect to
the standard Riemannian metric $g_1$
with exactly one zero in $x.$ 
The flow lines of $v_1$ consist of the circles passing 
through $x$ with a common tangent vector, see the pictures. 
Hence we obtain the following properties of the flow
$\phi^t: S^n \rightarrow S^n$ 
of the conformal vector field $v_1$ on $S^n:$ \\

 \begin{figure} 
 \centerline{{{\psfig{figure=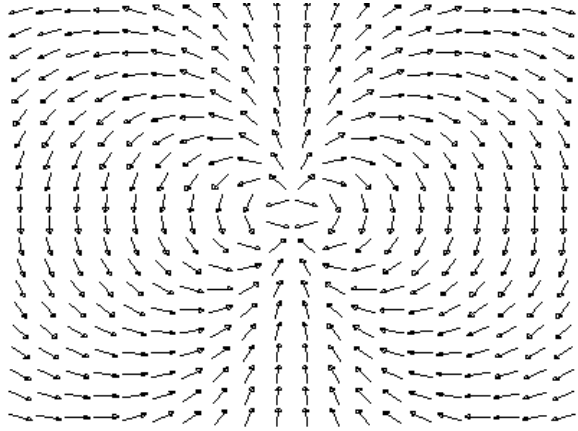,scale=0.4}}}}  
\caption{ The  vector field  ${v_1}/{|v_1|}$ in dimension 2  }
\end{figure}
\begin{figure} 
\centerline{{\psfig{figure=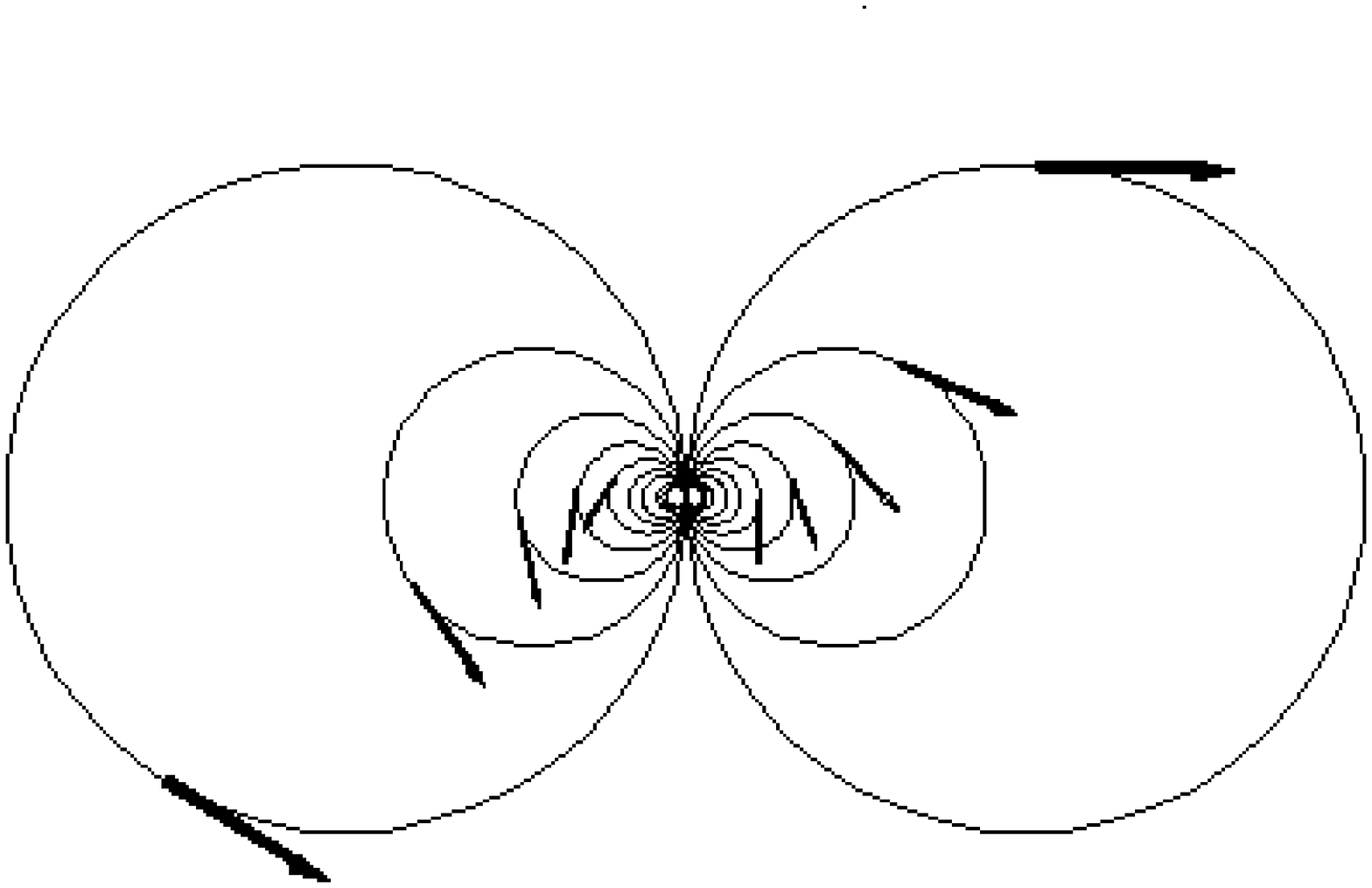,scale=0.4}}}  
\caption{   The integral curves of $v_1$ in dimension 2  }   
 \end{figure}

\begin{remark}
\label{rem:flow}
The flow $\phi^t$ of the conformal vector field $v_1$ defined
above satisfies the following properties:
\begin{itemize} 
\item[(a)] For any point 
$q \in S^n: x = \lim_{t \to \pm \infty} \phi^t(q).$
\item[(b)] For any tangent vector $\xi \in T_xS^n, F(x,\xi)=1$ 
there is a sequence
$q_i \in S^n-\{x\}$ with $\lim_{i\to \infty} q_i=x$ and
$\xi = \lim_{i\to \infty} \frac{v_1(q_i)}{F(q_i,v(q_i))}$ 
\end{itemize}
\end{remark}
Now we show that the mapping $h_A$ is an isometry also
for the Finsler metric $F_-:$ 
We can choose a sequence $m_i \to \infty$ with
$A^{m_i} \to 1.$ In particular for
a given $(z,\xi) \in T_z\rn; z\not= 0$ there is a unique 
$(0,\xi_0)\in T_0\rn; F((0,\xi_1))=1$ such that 
\begin{equation}
\label{eq:xinull}
\left(0,\xi_1\right) =
\lim_{i\to \infty} 
\frac{d\,\overline{f}_{A,b}^{m_i}(z,\xi)}{F_{-}
\left(d\,\overline{f}_{A,b}^{m_i}(z,\xi)\right)}\,.
\end{equation}
Since the mapping $\overline{f}_{A,b}$ is conformal 
for $F_-$ and since $h_A$ and $f_{A,b}$ commute it follows
that
\begin{eqnarray*}
\frac{F_{-}\left(dh_A (z,\xi)\right)}{F_{-}\left((z,\xi)\right)}
=\lim_{i\to \infty} 
\frac{F_{-}\left( 
d\,\overline{f}_{A,b}^{m_i}\,dh_A\,(z,\xi )\right)}
{F_{-}\left(d\,\overline{f}_{A,b}^{m_i}(z,\xi)\right)}
=\\
\lim_{i\to \infty} 
\frac{F_{-}\left(dh_A 
\left(d\,\overline{f}_{A,b}^{m_i}(z,\xi )\right)\right)}
{F_{-}\left(d\,\overline{f}_{A,b}^{m_i}(z,\xi)\right)}
= 
\frac{F_{-}\left(d\,h_A \left(0,\xi_0\right)\right)}{
F_{-}\left( \left(0,\xi_0\right)\right)}
=
\frac{F_{-}\left(d\,\left(\overline{f}_{A,b}\right)_0 
\left(0,\xi_0\right)\right)}{
F_{-}\left( \left(0,\xi_0\right)\right)}
=1
\end{eqnarray*}
as shown above, cf.
Equation~\ref{eq:fnull}.
Therefore the mapping $h_A$ is an isometry of 
the Finsler metric $F_{-}$. This implies that also the
flow generated by 
$\overline{f}_{1,b}=\overline{f}_{A,b} \circ h_A^{-1}$ 
is conformal for the Finsler metric
$F_{-}.$ Therefore the vector field $v_1$
on $S^n$ is also a conformal vector field for the 
Finsler metric $F$ on $S^n.$ 

Let us now  consider the following functions 
$m, M: S^n \rightarrow \R_{\ge 0}$: 
$$ 
m(q):= \frac{F^2(q, v_1(q))}{g_{(q)}(v_1(q), v_1(q))}, 
\  \ M(q):=\max_{\eta\in T_qS^n, \ \eta \ne 0}  
\frac{F^2(q, \eta)}{g_{(q)}(\eta, \eta)}- 
\min_{\eta\in T_qS^n, \ \eta \ne 0} 
\frac{F^2(q, \eta)}{g_{(q)}(\eta, \eta)} .  $$
Both functions are continuous functions invariant with respect 
to the flow $\phi^t$ of $v_1$. 
It follows from  Remark~\ref{rem:flow}(a) that
the function $m$ is a constant, i.e., there exists 
$\mu>0$ such that  
$F^2(q,v_1(q))=\mu \, g_{(q)}(v_1(q),v_1(q)).$
Part (b) of Remark~\ref{rem:flow} implies that  for every  
$0\ne \eta \in T_{x}S^n$ we have 
$\frac{F^2(x,\eta)}{g_{(x)}(\eta, \eta)}=\mu.  $
 Hence,  
$$M(x)= \max_{\eta\in T_qS^n, \ \eta \ne 0}  
\frac{F^2(q, \eta)}{g_{(q)}(\eta, \eta)}- 
\min_{\eta\in T_qS^n, \ \eta \ne 0} 
\frac{F^2(q, \eta)}{g_{(q)}(\eta, \eta)}=\mu- \mu =0.$$ 
But since $M$ is also flow invariant
by  Remark~\ref{rem:flow}(a), we have  
$M (q)=0$ for all $q \in S^n,$ i.e.,  
$F$ is up to a constant the norm of the
standard metric $g.$ 
Theorem 1 is proved.  
\\[1ex] 
\qed 

 \bigskip

As a consequence of the Proof of 
Theorem~\ref{thm:main} the inversion of the averaged Riemannian
metric is not a conformal map for a non-Euclidean Minkowski
metric, cf. Section~\ref{3}. 
Therefore one obtains from Liouville's theorem on
the conformal transformations of an Euclidean vector space
the following description of the
conformal transformations of a Minkowski space:
\begin{remark}
\rm
\label{rem:finsler-liouville}
Let $V$ be an $n$-dimensional vector space with a 
Minkowski norm $F$ which is not
Euclidean. Denote by $g$ the corresponding
averaged Euclidean metric. If $f:(U,F) \rightarrow (V,F)$ 
is a conformal mapping from an open
subset $U$ and $n\ge 3$ then $f$ is a   a similarity 
 with respect to the Minkowski metric
$F$ and with respect to the Euclidean metric $g$. 
Hence it is of the form
$ x\in V \mapsto \mu A x +b \in V$ for some
$\mu >0; b \in V$ and an orthogonal mapping $A$
of $(V,g).$
\end{remark}

\section{Conclusion}

Theorem 1 describes complete conformal vector fields 
of  Finsler metrics; it appears that no new  phenomena 
(with respect to  the Riemannian case) appear. 
Our proof is   based on the construction of averaged metric in Section \ref{sec2}, and 
 on  the description of  conformal vector 
fields for Riemannian metrics due to  
\cite{Lie,Li, Ob,Al,Yo,La,FT,Fe2,schoen}.

Let us also note that the existence of  a conformal vector 
field such that, for a certain 
point $p_0$, the closure of every trajectory contains 
this point is not artificial: as we know now, 
in view of Theorem 1, it is always the case, if the 
conformal transformations are essential.  
For a closed manifold, one also can show  it directly 
by repeating the Riemannian proof of \cite{Al}. 
  
As an interesting and much more involved problem in  
Finsler geometry related to transformation groups 
we would like to suggest to generalize the 
projective Lichnerowicz-Obata  conjecture for Finsler metrics, 
see \cite{Ma1,Ma2} for the proof of the Riemannian version, see also \cite{Sh}. \label{conc}
\\[2ex]
{\em Acknowledgement:} 
This work began when the authors met at the 
{\em 80\`eme rencontre entre physiciens th\'eoriciens et 
math\'emathiciens: G\'eom\'etrie de Finsler}
at Strasbourg University in September 2007.
We are grateful in particular to Athanase Papadopoulos for
the invitation to this conference. We also thank M. Eastwood, 
D. Hug,  Z. Shen and R. Schneider for useful discussions. 
We are also grateful to the referee for his suggestions.
The first  two authors thank   
Deutsche Forschungsgemeinschaft
(Priority Program 1154 --- Global Differential Geometry) for partial financial support.


{\small

Vladimir S. Matveev, Mathematisches Institut, Friedrich-Schiller Universit\"at Jena\\
07737 Jena, Germany, {\tt matveev@minet.uni-jena.de}

Hans-Bert Rademacher, Mathematisches Institut, Universit\"at Leipzig, 04081 Leipzig,
Germany \\
{\tt rademacher@math.uni-leipzig.de}

Marc Troyanov, Section de Math\'ematiques, {\'E}cole Polytechnique F{\'e}derale de
Lausanne\\
1015 Lausanne, Switzerland, {\tt marc.troyanov@epfl.ch}

Abdelghani Zeghib, UMPA, ENS-Lyon, 46, all{\'e}e d'Italie, 69364 Lyon Cedex 07, France\\
{\tt zeghib@umpa.ens-lyon.fr}

}
\end{document}